\documentstyle[12pt]{article}
\newcommand{\const}{\mathop{\rm const}\limits}

\newcommand{\Law}{\mathop{\rm Law}\limits}

\newcommand{\supp}{\mathop{\rm supp}\limits}

\begin{document}

\begin{center}

\vspace{3mm}

{\bf NON-IMPROVED UNIFORM TAIL ESTIMATES } \\

\vspace{3mm}

{\bf FOR NORMED SUMS OF INDEPENDENT RANDOM}\\

\vspace{2mm}

{\bf  VARIABLES WITH HEAVY TAILS,} \\

\vspace{2mm}

{\bf with applications.}\\

\vspace{3mm}

   E.Ostrovsky and L.Sirota, {\sc Israel}. \\

\vspace{2mm}

{\it Department of Mathematics and Statistics, Bar-Ilan University,
59200, Ramat Gan, Israel.}\\
e \ - \ mails: galo@list.ru; \ sirota@zahav.net.il \\

\vspace{3mm}
                     {\sc Abstract.}\\

\end{center}
\vspace{2mm}

 We obtain an uniform tail estimates for natural normed sums of independent random
 variables (r.v.) with regular varying tails of distributions.\par
  We give also many examples on order to show the exactness of offered estimates
and discuss some applications in the method Monte-Carlo and statistics, and obtain
the sufficient conditions for Central and stable limit theorem in the Banach space
of continuous function. \par
 There are considered a slight generalization on a random variables with super-heavy
tails and  martingale difference scheme.\par
\vspace{2mm}

\par

\vspace{2mm}

\par

\vspace{3mm}

\section{Introduction. Notations. Statement of problem.}
\vspace{3mm}

 Let $  \xi = \xi_1 $ be a random variable  defined on some sufficiently rich
 probabilistic space $ (\Omega, \cal{A}, {\bf P}) $ with regularly varying as $ x \to \infty $
 tail behavior:

 $$
 T(x) = T_{\xi}(x) = x^{-r} \ \log^{\gamma}(x) \ L(\log x),  r = \const > 0,
 \ x > e, \eqno(1.1)
 $$
where as ordinary the tail function $ T(x) = T_{\xi}(x)$ for the r.v. $  \xi $ may be
defined as follows:

$$
T(x) = T_{\xi}(x) = {\bf P}(|\xi| \ge x), \ x > 0; \eqno(1.2)
$$
and $ L = L(x) $ is positive continuous slowly varying as $ x \to \infty $ function. \par
 Further, we denote as $ \phi(t)= \phi_{\xi}(t) $ the characteristical function of the
 r.v. $\xi: $

$$
\phi(t)= \phi_{\xi}(t) = {\bf E} \exp(it\xi) \eqno(1.3)
$$
and by $  \psi(t) = \psi_{\xi}(t) $ its addition:

$$
\psi(t) = \psi_{\xi}(t) = 1 - {\bf Re}[ \phi_{\xi}(t)] =
{\bf E}(1- {\bf Re}[ \exp(it\xi)]). \eqno(1.4)
$$
 Note that if the r.v. $ \xi $ has symmetrical distribution, then

$$
\psi(t) = \psi_{\xi}(t) = 1 -  \phi_{\xi}(t) =
{\bf E}(1-  \exp(it\xi)). \eqno(1.4a)
$$

 Denote also in symmetrical case

 $$
 \overline{\psi}(t) = \overline{\psi}_{\xi}(t) =
 \sup_{\lambda \in (0,1)} \left[ \frac{\psi(\lambda t)}{\psi(\lambda)} \right];
 $$

 $$
 K(p) = 2 \pi^{-1} \Gamma(1+p) \sin(\pi p/2), \ p \in (0,2).
 $$

 Note that the function
 $$
  d(t,s) = d_{\xi}(t,s) = \sqrt{ \psi(t-s) }
 $$
 is a bounded translation invariant continuous {\it distance} between a two points on
 the real line  $  R. $ \par
  Let $ \xi(k), \ k=1,2,\ldots,n $ be independent copies of $ \xi; $
  we define the {\it natural norming} sequence $ \{ b(n) \} $ as follows: $ b(1) = 1 $
 and for $ n=2,3,\ldots $  as a positive solution of an equation:

$$
\psi_{\xi} (1/b(n)) = n.\eqno(1.5)
$$
 M.Braverman in \cite{Braverman2} proved that the non-random sequence $  \{b(n)\} $ is
 natural norming sequence for the sum of independent variables $ \{ \xi(i) \}. $ \par

 Note that for sufficiently greatest values $ n $ the values $ \{  b(n) \} $ exists,
 is unique  and $ \lim b(n) = \infty, \ n \to \infty. $ \par
  For instance, if $ {\bf Var}(\xi) \in (0,\infty), $ then  $ b(n) \asymp \sqrt{n}, $
 (the classical norming sequence.) \par
 Another example. Assume that the r.v. $  \xi $ has a symmetric stable distribution of
 order $ r: $
 $$
 \Law(\xi) = St(r) \ \Leftrightarrow \phi_{\xi}(t) =
 \exp \left( - |t|^r  \right), \ r \in (0,2).
 $$
 In this case $ b(n) \asymp \ n^{1/r}. $ \par

 Let us denote

$$
S(n) = b(n)^{-1} \sum_{k=1}^n \xi(k), \eqno(1.6)
$$

$$
U_{\xi}(x) = U(x) = \sup_{n=1,2,\ldots} T_{S(n)}(x). \eqno(1.7)
$$
 The function $ U(x) $ is {\it uniform tail function}  for natural normed sums of
 independent random variables $ \xi(k), \ k=1,2,\ldots,n. $\par
 \vspace{3mm}
{\bf Our aim is estimate of the uniform tail function $ U(x) $  through
the source tail function $  T(x). $ } \par
\vspace{3mm}
 We will distinguish and investigate all the cases $ r \in (0,2) $ (heavy tails),
 $  r = 2 $ (intermediate tails), and $ r > 2 $ (moderate tails);
  $ {\bf E}|\xi|^r = \infty  $ (infiniteness of main moment)
 and  $ {\bf E}|\xi|^r < \infty $ (finiteness of main moment). \par

 In the last section we consider the case of super-heavy tails, i.e. when
 the tail function decreases logarithmicaly.\par

  For the {\it exponential } tail function
 $$
 T(x) = \exp \left( -C x^m \right),  \ m=\const \in [1, \infty)
 $$
the non-improved estimate of the function $ U(x) $ has a view:

$$
U(x) \le \exp \left( -C_1(C,m) x^{\min(m,2)} \right),
$$
see  \cite{Braverman1}, \cite{Kozatchenko1}, \cite{Ostrovsky1}.\par
 Notice that in exponential case the critical value of the parameter $ m $ is also
 the value $ m=2.$ \par
  The {\it moment} estimates for $ S(n) $ in the case when
$ r \le 2 $ and $ {\bf E}|\xi|^r < \infty  $ is investigated in \cite{Braverman2}.
 Another estimates (tail and moments) see in \cite{Bahr1}, \cite{Carothers1},
\cite{Hitczenko1}, \cite{Johnson1}, \cite{Kwapien1}, \cite{Litvak1}, \cite{Prohorov1},
\cite{Vakhania1} etc. \par

 The applications of these estimates, e.g. in the Monte-Carlo method
 for errors estimates for integrals with infinite dispersion see in \cite{Egishjanz1},
\cite{Ostrovsky1}. In detail, let us consider the problem of numerical computation
of an absolute convergent integral (multiple, in general case) of a view:

$$
I = \int_D f(y) \ \nu(dy),
$$
where $ \nu(\cdot) $ is probabilistic measure on the set $ D: \ \nu(D) = 1.$ \par
 Let $ \tau(k), \ k=1,2,\ldots,n $ be  independent r.v. with distribution
 $ \nu: \ {\bf P}(\tau(k) \in A) = \nu(A). $  The Monte-Carlo estimation $ I_n $
 of an integral  $ I $ is

 $$
 I_n = n^{-1} \sum_{k=1}^n f(\tau(k)).
 $$

 Suppose for some $ r \in (1,2)$
$$
 {\bf E}|f(\tau(1))|^r  < \infty
$$
or more generally that the r.v. $ f(\tau(k)) - I $ satisfies the condition (1.1); note that
in the case $ r > 2 $ it may be used for error evaluating the classical Central Limit Theorem.\par
 In order to construct a non-asymptotical confidence interval for $ I $ of a reliability
 $ 1 - \delta, \ \delta = 0.05; 0.01 $ etc. we consider the probability

 $$
 U_n(x) = {\bf P} \left( b(n)^{-1} \left|\sum_{k=1}^n (f(\tau(k)-I)) \right| > x \right).
 $$
Note that $ U_n(x) \le U(x); $ therefore, we conclude denoting by $ X(\delta) $ the solution of
an equation

$$
U(X(\delta)) = \delta
$$
that with probability at least $ 1-\delta $

$$
|I_n - I| \le X(\delta) b(n)/n.
$$
 Obviously, under condition (1.1)
 $$
  \lim_{n \to \infty} b(n)/n = 0.
 $$

\vspace{3mm}

 Analogous application appears in statistics. Indeed, let us consider the following
classical  scheme of date elaboration.

$$
\eta(k) = \theta + \xi(k), \ k=1,2,\ldots,n;
$$
where $ \theta $ is unknown deterministic parameter, $ \{ \xi(k) \} $ are i.,i.d.
centered r.v. satisfying the condition (1.1) with $ r > 1 $ (additive noise with heavy
tail). \par
 The consistent estimation of the parameter $ \theta $ has a view

 $$
\hat{\theta}_n = n^{-1} \sum_{k=1}^n \eta(k).
 $$
We conclude as before that with probability at least $ 1-\delta $
under formulated above conditions and notations 
$$
\left|\hat{\theta}_n - \theta \right| \le X(\delta) b(n)/n.
$$

\vspace{3mm}

{\bf Remark 1.1.} Suppose

$$
T(x) = 0\left(x^{-r} \right), \ \ x \to \infty,  \ r = \const \in (0,2),
$$
then the distribution of $ \xi(k) $ dominated by symmetric stable distribution
$ C \cdot St(r); $ then by virtue of a main result of S.Kwapien \cite{Kwapien1}
 $ b(n) \sim n^{-1/r} $ and

$$
U(x) = 0\left(x^{-r} \right), \ r = \const \in (0,2).
$$

{\bf Remark 1.2.} Suppose

$$
{\bf E} |\xi|^r < \infty, \ r = \const \in (0,2).
$$
 This condition is equivalent to the follows:

 $$
 \int_e^{\infty} x^{-1} \ \log^{\gamma}(x) \ L(\log x) \ dx < \infty.
 $$

  We conclude using famous result belonging to B. von Bahr and C.-G. Esseen \cite{Bahr1}
that again $ b(n) \sim n^{-1/r} $ and

$$
\sup_n {\bf E}|S(n)|^r \le 2 {\bf E} |\xi|^r,
$$
therefore

$$
U(x) = o\left(x^{-r} \right), \ x \to \infty, \ r = \const \in (0,2).
$$

\vspace{3mm}
  The paper is organized as follows. In the next section we consider the so-called case
 of heavy tails: $  r \in (0,2). $  Third section is devoted mainly to the consideration of
 intermediate case $ r=2. $ Fourth section contains the investigation of the case moderate
 tails: $ r \in (2,\infty). $ \par
  In the fifth section we generalize preceding results on the martingale case, i.e. when
  the summands $ \{\xi(i) \} $ are centered martingale differences relative some filtration
  $  \{ F(i) \}.$  The sixth section  contains the tail evaluating for normed sums of the r.v.
  with superheavy, i.e. logarithmical, tails of distribution. In the $ 7^{th} $     section
  we obtain sufficient conditions for  Stable and  Central Limit Theorems for heavy tail
  random fields in the Banach space of continuous functions on the compact metric spaces. \par
    The last section is devoted to concluding remarks.

\vspace{3mm}

The letter $ C, $ with or without subscript, denotes a finite positive non essential constants,
not necessarily  the same at each appearance. \par

\vspace{3mm}

\section{Main result: heavy tails}

\vspace{3mm}

{\bf Theorem 2.1.} Let the r.v. $ \xi $ has a symmetrical distribution. Then

$$
U(x) \le  0.5  x \int_{-2/x}^{2/x} \overline{\psi}(t) \ dt, \ x > 0. \eqno(2.1)
$$
{\bf Proof.}  We will use the well-known inequality, which is true for arbitrary
r.v. $ \eta: $

$$
T_{\eta}(2/a) \le a^{-1}\int_{-a}^a \psi_{\eta}(t) \ dt, \ a = \const > 0. \eqno(2.2)
$$
 Here $ a = 2/x. $ We denote

 $$
 V(n) = \sum_{k=1}^n \xi(k) = b(n) S(n).
 $$
 M.Braverman in  \cite{Braverman2} proved that

 $$
 \psi_{V(n)}(t) \le \sum_{k=1}^n \psi_{\xi(k)}(t) = n \ \psi_{\xi}(t) = n \ \psi(t),
 $$
 therefore

 $$
 \psi_{S(n)}(t) \le n \ \psi(t/b(n)).
 $$
 We get:

 $$
 \psi_{S(n)}(t) \le \sup_n [n \ \psi(t/b(n))] \le \sup_n [n \ \psi(t \cdot \psi^{-1}(1/n))] \le
 $$

 $$
 \sup_{\lambda \in (0,1)} \frac{\psi(\lambda t)}{\psi(\lambda)} = \overline{\psi}(t). \eqno(2.3)
 $$

 It remains to use the inequality (2.2). \par

\vspace{3mm}
 {\bf Theorem 2.2.} Suppose again the r.v. $ \xi $ has a symmetrical distribution and let
 $ r \in (0,2).$ Then
$$
U(x) \le \inf_{p \in (0,r)} \left[K(p) \ x^{-p} \
\int_0^{\infty} \overline{\psi}(t) \ t^{-p-1} \ dt \right]. \eqno(2.4)
$$
{\bf Proof.}  We will use the famous formula:

$$
{\bf E} |\eta|^p = K(p) \int_0^{\infty} \psi_{\eta}(t) \ t^{-p-1} \ dt, \eqno(2.5)
$$
see \cite{Braverman2}, \cite{Kawata1}. Here $ \eta = S(n) $ and we may replace in the last
inequality $ \overline{\psi}_{S(n)}(t) $   instead $ \psi_{\eta}(t).$\par
 Estimating $ \overline{\psi}_{S(n)}(t) $  by means of inequality (2.3), we obtain the
 assertion of theorem 2.2.\par

{\bf Corollary 2.1. } Choosing in (2.4) the value $ p = r-C/\log x, \ x > e^2, $ we conclude:

$$
U(x) \le e^{C} \left[K(r-C/\log x) \ x^{-r} \
\int_0^{\infty} \overline{\psi}(t) \ t^{-r+C/\log x-1} \ dt \right]. \eqno(2.6)
$$

 If for example
$$
 \overline{\psi}(t) \le |t|^r \ |\log t|^{\beta}, \ |t| < 1/e, \ \beta = \const \ge 0, \eqno(2.7)
$$
then the optimal value of the constant $ C $ in (2.5) is $ C=\beta+1 $ and hence
for the values $ x  > \exp(2(\beta+1)/r) $

$$
U(x) \le e^{\beta+1} \ (\beta+1)^{-\beta-1} \ K(r-(\beta+1)/\log x) \ x^{-r} \
[\log x]^{\beta+1} \le
$$

$$
 C_1(\beta,r) \ x^{-r} \ [\log x]^{\beta+1}.\eqno(2.8)
$$

\vspace{3mm}

{\it We investigate further in this section the case only when in the condition (1.1)
$ r \in (0,2). $ } \par

\vspace{3mm}

{\bf Definition 2.1.} We will say that the r.v. $ \eta $ has a {\it regular tail,} write:
$ \Law(\eta) \in RT, $ if is true
the inverse inequality to the inequality (2.2) up to multiplicative constant, i.e. when there
exists constant $ C_1 > 0 $ so that for all sufficiently small positive values $  a $

$$
T_{\eta}(2/a) \ge C_1 a^{-1}\int_{-a}^a \psi_{\eta}(t) \ dt. \eqno(2.9)
$$
 For instance, if the non-trivial r.v. $ \xi $ satisfies the condition (1.1) with $ r \in (0,2], $
 then $ \Law(\xi) \in RT. $ \par
  Conversely, if $ \xi \ne 0 $ and if for some $ \Delta = \const > 0 \ {\bf E}|\xi|^{2+\Delta} < \infty, $
  then $ \Law(\xi) \notin RT. $ Namely, it follows from Tchebychev's inequality

$$
T_{\xi}(2/a) \le C_2 a^{2+\Delta}, \ a \in (0, 1),
$$
 but

$$
 a^{-1}\int_{-a}^a \psi_{\xi}(t) \ dt \ge C_3 a^2.
$$

\vspace{3mm}

{\bf Definitions 2.2.} We will say that the r.v. $ \xi $ (more exactly, the distribution $ \Law(\xi) $
of the r.v. $ \xi) $ belongs to the class $ MI $ (Monotonically Increasing), write: $ \Law(\xi) \in MI, $
if  the function $ \lambda \to \theta(\lambda), \ \lambda \in (0,1), $ where

$$
\theta(\lambda)= \frac{\psi(\lambda t)}{\psi(\lambda)}, \eqno(2.10)
$$
monotonically non-decreased for all the values $  t  $ in some neighborhood
$ t \in [0, \Delta), \ \Delta > 0. $ \par

 Analogously, we will say that the r.v. $ \xi $ (more exactly, the distribution $ \Law(\xi) $
of the r.v. $ \xi) $ belongs to the class $ MD $ (Monotonically Decreasing), write: $ \Law(\xi) \in MD, $
if  the function $ \lambda \to \theta(\lambda), \ \lambda \in (0,1), $ where

$$
\theta(\lambda)= \frac{\psi(\lambda t)}{\psi(\lambda)},\eqno(2.11)
$$
monotonically non-increased for all the values $  t  $ in some neighborhood
$ t \in [0, \Delta), \ \Delta > 0. $ \par

 Obviously, if $ \Law(\xi) \in MD, $ then

 $$
 \overline{\psi}(t)= \psi(t), \ t \in [0,\Delta] \eqno(2.12)
 $$
 and if $ \Law(\xi) \in MI $ and if the distribution of the r.v. $\xi$ satisfies the condition (1.1), then

 $$
 \overline{\psi}(t)= t^r,   \ t \in [0,\Delta]. \eqno(2.13)
 $$

 For instance,
 if the distribution of the r.v. $\xi$ satisfies the condition (1.1)
 and $ \gamma > 0, $ then $ \Law(\xi) \in MD $ and $ \overline{\psi}(t)= \psi(t), \ t \in [0,\Delta];$
if the distribution of the r.v. $\xi$ satisfies the condition (1.1)
 and $ \gamma < 0, $ then $ \Law(\xi) \in MI $ and $ \overline{\psi}(t)= t^r, \ t \in [0,\Delta].$\par

 Consider now the case $ \gamma = 0. $ Suppose

 $$
 T(x) = T_{\xi}(x) = x^{-r} \ [\log \log x]^{\kappa} \ L(\log\log x),  \kappa = \const,
 \ x > e^3, \eqno(2.14)
 $$
where as before $ L = L(z) $ is slowly varying as $ z \to \infty $ positive function. \par
 If in (2.13) $ \kappa >0, $ then $ \Law(\xi) \in MD $ and if
 in (2.13) $ \kappa < 0, $ then $ \Law(\xi) \in MI. $\par

\vspace{3mm}

{\bf Theorem 2.3.} If $ \Law(\xi) \in RT \cap MD, $ then we propose the following
non-improvable up to multiplicative constant estimates:

$$
T(x) \le U(x) \le C(r,L) T(x), \ x > x_0 = \const > 1. \eqno(2.15)
$$

\vspace{3mm}

{\bf Proof.} The left inequality in (2.14) is trivial; it remains to prove the right-hand inequality.\par
We will use the following fact: if the equality (1.1) holds and $ r \in (0,2), $ then as $ t \to 0+ $

$$
\psi_{\xi}(t) \sim C_3(r) \ t^r \ |\log t|^{\gamma} \ L(|\log t|), \eqno(2.16)
$$
where

$$
C_3(r) = \Gamma(1-r) \ \cos(\pi r/2), \ r \ne 1,
$$
\cite{Ibragimov1}, p. 86-87, see also \cite{Braverman2}.\par
  As long as $ \Law(\xi) \in MD, $

 $$
 \overline{\psi}_{\xi}(t) \sim {\psi}_{\xi}(t), \ t \in (0,\Delta).\eqno(2.17)
 $$
Since $ 0 < r < 2, $ we conclude by virtue of the equality (2.16) and the condition  $ \Law(\xi) \in RT $
that

$$
U(x) \le C(r,L) x^{-r} \ [\log x]^{\gamma} \ L(\log x) = C(r,L) T(x),\eqno(2.18)
$$
Q.E.D.\par

 Analogously may be proved the following result.

\vspace{3mm}

{\bf Theorem 2.4.} If $ \Law(\xi) \in RT \cap MI, $ then

$$
 U(x) \le C_1(r,L) x^{-r}, \ x > x_0 = \const > 1, \eqno(2.19)
$$
and the last inequality is exact, for example, for the symmetric stable distribution $ \xi. $\par

 Notice that the member $ L(\cdot) $ is absent in the right-hand of the inequality (2.19).
 This expression included in the definition on the norming sequence $ \{b(n)\}. $\par

\vspace{3mm}

\begin{center}
{\bf Examples.}\par
\end{center}
{\bf A.} Suppose the r.v. $  \xi $ is symmetrically distributed,
 satisfies the condition (1.1) and $  \gamma > 0; $ then

 $$
 U(x) \le  C(r,\gamma,L) x^{-r} \ [\log x]^{\gamma} \ L(\log x) = C(r,L) T(x), \ x > e^2.
 $$

{\bf B.}  Suppose now the r.v. $  \xi $ is symmetrically distributed,
 satisfies the condition (1.1) and $  \gamma < 0; $ then

$$
 U(x) \le C_1(r,L) x^{-r}, \ x > x_0 = \const > 1.
$$

{\bf C.} Consider now up to end of this section the case when $ \gamma = 0. $ Suppose

 $$
 T(x) = T_{\xi}(x) = x^{-r} \ [\log \log x]^{\kappa} \ L(\log\log x),  \kappa = \const,
 \ x > e^3,
 $$
where as before $ L = L(z) $ is slowly varying as $ z \to \infty $ positive function. \par
 If  $ \kappa >0, $  then

 $$
 U(x) \le C(r,\kappa,L) x^{-r} \ [\log \log x]^{\kappa} \ L(\log\log x), \ x > e^3.
 $$

{\bf D.} If  $ \kappa < 0, $ we conclude

 $$
 U(x) \le C_1(r,\kappa,L) x^{-r}, \ x > 1.
 $$

 {\bf Remark 2.1.} Note that the application of theorem 2.2. give us more
 slight result, namely

 $$
 U(x) \le C \ x^{-r} \ [\log x]^{\gamma+1} \ L(\log x) = T(x), \ x > e.
 $$

\vspace{3mm}

\section{Main result: intermediate case.}

\vspace{3mm}
{\it We consider in this section the case when in the condition (1.1) $ r=2. $ } \par

\vspace{3mm}
{\bf Theorem 3.1.} Suppose the r.v. $  \xi $ is symmetrically distributed,
 satisfies the condition (1.1) for $ r=2 $ and for some $  \gamma =\const \in R; $ then: \par
 {\bf A.} If $ \gamma \ge -1, $ then

 $$
 U(x) \le  C(L) x^{-2} \ [\log x]^{\gamma+1} \ L(\log x) = C(r,\gamma,L) \ T(x) \ \log x, \ x > e; \eqno(3.1a)
 $$

 {\bf B.} If $ \gamma < -1, $ then

 $$
 U(x) \le  C(L) x^{-2}, \ x > 1. \eqno(3.1b)
 $$

{\bf Proof.} Let us introduce a following function:

$$
H(x) = - \int_0^x u^2 \ d T_{\xi}(u); \eqno(3.2)
$$
then
$$
\psi_{\xi}(t) \sim 0.5 \ t^2 \ H(1/|t|), \ t \to 0,
$$
see \cite{Ibragimov1}, p.86-88; see also \cite{Braverman2}.\par
 In the considered case

$$
\psi_{\xi}(t) \sim C_4(\gamma,L) \ t^2 \ |\log t|^{\gamma+1} \ L(\log x) , \ t \to 0. \eqno(3.3)
$$
 The remained part of the proof is at the same as in theorems 2.3-2.4.\par
  {\bf Remark 3.1.} Note concerning the lower bound for $ U(x) $  in considered case $ r = 2 $ that
we can show only the trivial bound

 $$
 U(x) \ge  x^{-2} \ [\log x]^{\gamma} \ L(\log x) = T(x), \ x > e; \eqno(3.4)
 $$
 Note that there is a "gap" of a view  "degree of a logarithmic term" $ [\log x]^\Delta, \ \Delta > 0 $
  between upper and lower bound for   the uniform tail of probability $ U(x). $\par

\vspace{3mm}

\section{Main result: moderate tails}

\vspace{3mm}

{\it We concentrate our attention in this section on the case when in the equality (1.1)
$ r > 2, $ and suppose $ {\bf E}\xi = 0. $} \par
 Recall that in this case the norming sequence $ b(n) $ is ordinary: $ b(n) = \sqrt{n}. $ \par
\vspace{3mm}

 In order to formulate and prove the main result in this case, we recall here for reader
 convenience some facts about so-called Grand Lebesgue Spaces \cite{Fiorenza1},
 \cite{Fiorenza2}, \cite{Iwaniec1}, \cite{Liflyand1} etc.
or equally  "moment"  spaces of random variables defined on
 fixed probabilistic space $ (\Omega, \cal{A}, {\bf P}); $ more detail description see
 in  \cite{Kozatchenko1}, \cite{Liflyand1}, \cite{Ostrovsky1}, \cite{Ostrovsky2}.\par

 Let us consider the following norm (the so-called "moment norm")
on the set of r.v. defined in our probability space by the following way: the space
$ G(\nu) = G(\nu;r) $ consist, by definition, on all the  r.v. with finite norm

$$
||\xi||G(\nu) \stackrel{def}{=} \sup_{p \in (2,r)} [|\xi|_p/\nu(p)], \ |\xi|_p :=
{\bf E}^{1/p} |\xi|^p.\eqno(4.1)
$$
  Here $ r = \const > 2, \ \nu(\cdot) $ is some continuous positive on the
 {\it semi-open} interval  $ [1,r) $ function such that

     $$
     \inf_{p \in (2,r)} \nu(p) > 0, \ \nu(p) = \infty, \ p > r;
     $$
 and as usually

 $$
 |\xi|_p \stackrel{def}{=}  \left[{\bf E} |\xi|^p \right]^{1/p}
 $$

We will denote
$$
 \supp (\nu) \stackrel{def}{=} [1,r) = \{p: \nu(p) < \infty \}.
$$

 The case $ r = +\infty $ is investigated in  \cite{Kozatchenko1}, \cite{Ostrovsky1},
 \cite{Ostrovsky2}; therefore, we suppose further $  2 < r < \infty. $\par

 Let $  \xi $ be a r.v. such that
 $$
 p > r \Rightarrow
 |\xi|_p \stackrel{def}{=} \left[ {\bf E}|\xi|^p \right]^{1/p} = \infty \
 $$
The {\it natural} function   $  \nu_{\xi}(p)  $ may be defined as follows:

$$
\nu_{\xi}(p) := |\xi|_p = \left[ {\bf E}|\xi|^p \right]^{1/p}.\eqno(4.2)
$$
 Obviously,

 $$
 ||\xi||G(\nu_{\xi})=1.
 $$

The {\it natural} function for the {\it family}  $ \{\xi(\cdot) \} $
of a r.v.  $ \{\nu(\cdot) \} = \{\nu_{\xi(k)}(p) \}, \ k = 1,2,\ldots  $ may be
defined as follows:

$$
\nu_{ \{\xi\}}(p) := |\xi|_p = \sup_k \left[ {\bf E}|\xi(k)|^p \right]^{1/p},\eqno(4.2a)
$$
if there exists  and is finite.

\vspace{4mm}

 The  complete description of a possible natural functions see in
\cite{Ostrovsky2}, \cite{Ostrovsky1}, chapter 1, section 3.\par

{\bf Example 4.1.} Suppose the r.v. $ \xi $ satisfies the condition (1.1) for $ r > 1 $ and
$ \gamma > -1. $ Then (see \cite{Liflyand1}, \cite{Ostrovsky2}) for  the values
$ p \in [1,r), \ p \to r-0 $

$$
{\bf E}|\xi|^p \sim C(r,\gamma,L) \ (r-p)^{-\gamma-1} \ L(1/(r-p)).
$$

  Note that an inequality $ \psi(r-0) < \infty $ is equivalent to the moment restriction
 $  |\xi|_r < \infty. $  \par

 We recall now the relations between moments for r.v. $  \xi, \ \xi \in G(\nu,r) $
 and its tail behavior. Namely, for $  p < r $

 $$
 |\xi|_p = \left[ p \int_0^{\infty} u^{p-1} \ T_{|\xi|}(u) \ du \right]^{1/p},
 $$
therefore

$$
||\xi||G(\nu;r) = \sup_{p < r}
\left\{\left[ p \int_0^{\infty} u^{p-1} \ T_{|\xi|}(u) \ du \right]^{1/p}/\nu(p) \right\}.
$$
 Conversely, if the r.v. $  \xi $ belongs to the space $ G(\nu,r), $ then

$$
T_{|\xi|}(x) \le T^{(||\xi||G(\nu;r) \cdot \nu)}(x),
$$
where by definition
$$
T^{(\nu)}(x) \stackrel{def}{=}
\inf_{p \in (1,r)} \left[ \nu^p(p)\ /x^p \right], \ x > 0.
$$

{\bf Example 4.2.} If

$$
{\bf E}|\xi|^p  \le C_2 \ (r-p)^{-\gamma-1} \ L(1/(r-p)), \ p < r,
$$
then

$$
T_{|\xi|}(x) \le C_4(r,\gamma,L) \ x^{-r} \ [\log x]^{\gamma+1} \ L(\log x), \ x > e.
$$
 Notice that there is  a "logarithmic gap" between  upper and lower tail and moment relations.
 This gap is essential, see \cite{Liflyand1}, \cite{Ostrovsky2}. Let us consider the following \par
 {\bf Example 4.3.} Let $ \zeta $ be a discrete r.v. with distribution

$$
{\bf P} \left( \zeta = \exp(e^k)  \right) = C_5 \ \exp \left(\beta r k - r  e^k \right), \
$$
$$
 k = 1,2,\ldots, \ \beta = \const > 0,
$$
where obviously

$$
1/C_5 = \sum_{k=1}^{\infty} \exp \left(\beta r k - r  e^k \right).
$$
 We conclude after some calculations:

 $$
 |\zeta|_p \asymp (r-p)^{-\beta},
 $$
but for the sequence $ x(k)= \exp(\exp(k)) $
$$
T_{\zeta}(x(k)) \ge c_6 [\log x(k)]^{\beta} \ (x(k))^{-r}.
$$

\vspace{3mm}
 Let us introduce  the so-called Rosenthal's constants (more exactly, Rosental's functions)
 \cite{Rosenthal1} $ R(p) $ as follows:

$$
R(p) = \sup_{\zeta(k): {\bf E}\zeta(k)=0}
\frac{|\sum_{k=1}^n \zeta(k)|_p}{\max(|\sum_{k=1}^n \zeta(k)|_2,
\left[ \sum_{k=1}^n |\zeta(k)|^p_p \right]^{1/p})},
\eqno(4.3)
$$
where the supremum in (4.3) is calculated over all the sequences of independent centered
random variables $ \{\zeta(k)\} $  with condition $ |\zeta(k)|_p < \infty. $ \par
 This constants was intensively investigated in many publications, see,e.g. \cite{Rosenthal1},
 \cite{Carothers1}, \cite{Hitczenko1},   \cite{Johnson1},  \cite{Johnson2}, \cite{Latala1},
 \cite{Ostrovsky3}, \cite{Sharachmedov1}, \cite{Utev1}  etc. It is known, for instance, that
 there exists an absolute constant $ C_R, $ which is exactly calculated in \cite{Ostrovsky3}:
$ C_R \approx  1.77638, $ so that

 $$
 R(p) \le C_R \frac{p}{e \cdot \log p}, \ p \in [2, \infty). \eqno(4.4)
 $$
Note that for symmetrical distributed r.v. $ \{\zeta(k) \} $ the correspondent Rosenthal's
 constant $ C_R$ is (approximately) equal to $ \approx 1.53573. $\par

 It follows from inequality (4.4) that if the r.v. $ \{ \eta(k)\}, \ k=1,2,\ldots,n  $ are
centered,  i., i.d., and $ |\eta(1)|_p < \infty, \ p \ge 2,  $ then

$$
\sup_n \left|n^{-1/2} \sum_{k=1}^n \eta(k) \right|_p \le R(p) \ |\eta(1)|_p. \eqno(4.5)
$$

 Suppose now that the {\it centered}  r.v. $  \xi $ belongs to some space
 $ G(\nu;r), \ r \in (2,\infty); $ this condition is satisfied, if $ \xi $  satisfied
 the condition (1.1) for some $ r \in (2,\infty)$ and $ {\bf E}\xi = 0. $ \par
 For instance, the function $ \nu = \nu(p) $ may be the natural function for the r.v.
 $ \xi: \nu = \nu_{\xi}(p). $ \par

 \vspace{3mm}

{\bf Theorem 4.1.} We have the following non-improved inequality in the terms of $ G(\nu;r) $
norms:

$$
 \sup_n ||n^{-1/2} \sum_{k=1}^n \xi(k) ||G(\nu;r)
 \le C(\nu;r) ||\xi||G(\nu;r); \eqno(4.6a)
$$

$$
||\xi||G(\nu_{\xi}) \le \sup_n ||n^{-1/2} \sum_{k=1}^n \xi(k) ||G(\nu_{\xi})
 \le C(\nu_{\xi}) ||\xi||G(\nu_{\xi}). \eqno(4.6b)
$$

\vspace{3mm}
{\bf Proof.} Let $ \xi \in G(\nu;r) $ and $ {\bf E}\xi = 0. $ We can and will assume
without loss of generality $ ||\xi||G(\nu;r)=1. $ It follows from this equality

$$
|\xi|_p \le \nu(p), \ p \in [2,r).
$$
 We conclude by means of inequality (4.5)  by virtue of inequality $ p < r: $

$$
 \sup_n \left|n^{-1/2} \sum_{k=1}^n \xi(k) \right|_p  \le R(p) |\xi|_p \le R(p) \
 \nu(p) \le C_4(r) \nu(p), \eqno(4.7)
$$

$$
\sup_n \sup_{p \in [2, r)}  \left[ \left|n^{-1/2} \sum_{k=1}^n \xi(k) \right|_p /\nu(p) \right] \le C_4(r),
$$
which is equivalent to the assertion (4.6a) of our theorem.\par
The proposition (4.6b) follows from (4.6a) after choosing $ \nu(p)= \nu_{\xi}(p). $\par

 As a consequence: \par

 \vspace{3mm}

 {\bf Theorem 4.2.}

 $$
 U(x) \le C(r,\Law(\xi)) \ T^{(\nu) }(x). \eqno(4.8)
 $$

\vspace{3mm}

{\bf Example 4.4.} Let $ {\bf E}\xi = 0 $ and suppose the r.v. $  \xi $ satisfies the condition (1.1).
Then

 $$
 x^{-r} \ \log^{\gamma}(x) \ L(\log x) \le U(x) \le
 C_6(r,\gamma,L) \ x^{-r} \ \log^{\gamma+1}(x) \ L(\log x), \ x > e. \eqno(4.9)
 $$

\vspace{3mm}
{\bf Remark 4.1.} {\it We can not prove that the power of the logarithmic multiplier} $ \gamma+1 $
{\it  in the right-hand side of  last inequality  is  not improvable. } \par
 But we can offer the following example.\par
 \vspace{3mm}
{\bf Example 4.5.} Let $ \zeta  $  be the r.v. from the example 4.3., put $ \zeta^o =
\zeta - {\bf E}\zeta. $ \par
Define a r.v. $ \theta $  as a r.v. with
{\it accomplishing} infinite divisible distribution
  to the distribution $ \zeta^o. $ This imply that

  $$
  \theta = \sum_{m=1}^{\tau} \zeta^o(m),
  $$
  where $  \zeta^o(m) $ are independent copies of  $ \zeta^o, \ $  the r.v. $ \tau $ gas a standard
  Poisson  distribution with unit parameter: $  {\bf P}(\tau=j) = 1/(e \ j!), j=0,1,2,\ldots; \
  \sum_{m=1}^0 = 0. $\par
   Obviously,

   $$
   \phi_{\theta}(t) = \exp \left[ \phi_{\zeta^o}(t) - 1 \right].
   $$
   It follows from the last formula that the asymptotical  behavior of the function
$ \phi_{\theta}(t)$ as $t \to 0  $ or equally the tail behavior $ T_{\theta}(x) $ as
$ x \to\infty $ alike the asymptotical  behavior of the function
$ \phi_{\zeta^o}(t), \ t \to 0 $ and correspondingly the behavior $ T_{\zeta^o}(x)$
as $ x \to \infty. $  For instance, Yu.V.Prohorov in  \cite{Prohorov1} proved the following
inequality:

$$
{\bf P} \left( \sum _{k=1}^n \theta(k) >x \right) \le 8 \
 {\bf P} \left( \sum_{k=1}^n \zeta^o(k) > x/2 \right). \eqno(4.10)
$$

\vspace{3mm}

 Let us denote

 $$
 \Theta_n = n^{-1/2}\sum_{k=1}^n \theta(k),
 $$
where $ \{ \theta(k) \}  $ are independent copies of the r.v. $ \theta. $  We conclude:

 $$
 |\theta|_p \asymp (r-p)^{-\beta}, \ p \in [1,r),
 $$
but there exists a constant $ q \in (0,1) $  such that
for the sequence $ x(k)= \exp(\exp(k)) $  and for all the values $ n=1,2,\ldots $
$$
T_{\Theta_n}(x(k)) \ge C_7 \ q^n \ [\log x(k)]^{\beta} \ (x(k))^{-r}.
$$

\vspace{3mm}
{\bf Remark 4.2.} It is not necessary to suppose in this section that the r.v. $ \xi(i) $ are
identical distributed; it is sufficient to assume that the r.v. $ \{\xi(i) \} $ are
independent,centered and  such that

$$
\sup_i T_{\xi(i)}(x) \le x^{-r} \ \log^{\gamma}(x) \ L(\log x),  r = \const > 2,
 \ x > e. \eqno(4.11)
 $$

\vspace{3mm}

{\bf Remark 4.3.} The application of the interpolation technique, see for instance
\cite{Karlowich1}, gives us more slight result as in the theorems  4.1-4.2. For example,
we conclude under conditions in example  4.4

 $$
  U(x) \le C_7(r,\gamma,L) \ x^{-r} \ \log^{\gamma+1}(x) \ \log \log x \
  L(\log x), \ x > e^e. \eqno(4.12)
 $$

 We give briefly here the proof of last inequality, as long as it has by our opinion
self-contained interest. \par
 Let the centered r.v. $ \xi(k), \ k=1,2,n $  be i., i.d. and satisfy the equality (1.1).
 Let also $ \delta = \const > 0. $ Consider the following Orlicz's function:

 $$
 N(u) = N_{r,\gamma, L}(u) = \delta^{-1} \ |u|^r \  \left( \log^{-\gamma - \delta} |u| \right) \
 L(\log |u|), \ u > e, \eqno(4.13)
 $$
and as usually $ N(u) = C \ u^2, |u| \le e, $ so that the function $ N = N(u) $ is continuous.\par
We can realize the sequence $ \xi(k), \ k=1,2, \ldots $ on the direct (Cartesian) product
of probabilistic spaces $ \left(\Omega_k, {\cal A}_k, {\bf P_k} \right): $

$$
(\Omega, {\cal A}, {\bf P}) = \otimes_{k=1}^{\infty} \left(\Omega_k, {\cal A}_k, {\bf P_k} \right)
$$
so that $ \xi(k) $ is symmetrical distributed measurable function $ \xi(k): \Omega_k \to R. $ \par

 As long as $ r > 2, $ there is two numbers $ p_1, p_2 $ for which $ 2 < p_1 < r < p_2 < \infty, $
for example, $ p_1 = 0.5 (2+r), \ p_2 = r+1. $ \par
 Define the following linear operator acting only on the centered random vectors:

 $$
 T(\xi(1), \xi(2), \ldots, \xi(n)) = n^{-1/2} \sum_{k=1}^n \xi(k). \eqno(4.14)
 $$
  From the Rosental's inequality follows that $  T $ is bounded operator from the spaces 
$ L_{p_j}, j=1,2 $ into  at the same spaces. We conclude by virtue of interpolation theorem 
belonging to A.Yu.Karlowich and L.Maligranda \cite{Karlowich1} that the operator $ T $
is bounded uniformly in $ n $ operator from the  Orlicz's space
$ M(N_{r,\gamma, L}, \Omega, {\bf P}) \stackrel{def}{=} M(r,\gamma, L) $
into at the same space:

$$
\sup_n ||S(n)||M(r,\gamma, L) \le C(r,\gamma, L) \ \sup_k ||\xi(k)||M(r,\gamma, L). \eqno(4.15)
$$
 It is easy to calculate that the right-hand side of last inequality is uniformly in
$ k $ and $ \delta \in (0, 0.08) $ bounded and therefore

$$
\sup_{\delta \in (0, 0.08)}  \sup_n ||S(n)||M(r,\gamma, L) \le C_2(r,\gamma, L) < \infty. \eqno(4.16)
$$

 Since the function $ N(u) = N_{r,\gamma, L}(u) $ satisfies the  $  \Delta_2  $ condition,
 we get: (see \cite{Krasnoselsky1}, chapter 2, section 9)

\vspace{3mm}

{\bf Proposition 4.1.}
$$
\sup_{\delta \in (0, 0.08)}  \sup_n {\bf E}N_{r,\gamma, L}(S(n)) \le C_2(r,\gamma, L)
\sup_{\delta \in (0, 0.08)} \sup_k {\bf E}N_{r,\gamma, L}(\xi(k)). \eqno(4.17)
$$

\vspace{3mm}
 Tchebychev's inequality give us:

$$
\sup_{\delta \in (0, 0.08)}  \sup_n T_{S(n)}(x) \le C_3(r,\gamma, L) \ x^{-r} \ \delta^{-1} \
\log^{\gamma+1 + \delta} x \ L(\log x).
$$
 We obtain choosing $ \delta = (\log \log x)^{-1} $

 $$
  U(x) \le C_7(r,\gamma,L) \ x^{-r} \ \log^{\gamma+1}(x) \ \log \log x \
  L(\log x), \ x > e^e,
 $$
Q.E.D.

\vspace{3mm}

\section{Martingale generalization}

 In this fifth section we generalize preceding results on the martingale case, i.e. when
  the summands $ \{\xi(i) \} $ are centered $ {\bf E}\xi(i) = 0 $ martingale differences
  relative some filtration $  \{ F(i) \}: \ F(0) = \{ \emptyset, \ \Omega \},  F(i) \subset
  F(i+1) \subset \cal{A}.$ This imply, by definition,

  $$
  {\bf E} \xi(k)/F(k-1)=0, \ k=1,2,\ldots.
  $$

  \vspace{3mm}

 We will use in the sequel the following generalization of Rosenthal's inequality
for centered martingales, see \cite{Ostrovsky4}:

$$
 \left|n^{-1/2}\sum_{k=1}^n \xi(k) \right|_p \le p \ \sqrt{2} \ \sup_i |\xi(i)|_p, \ p \ge 1.  \eqno(5.1)
$$

{\bf Theorem 5.1.} We have alike in the fourth section in the case $ 2 < r < \infty $
the following non-improved inequality in the terms of $ G(\nu;r) $ norms:

$$
 \sup_n ||n^{-1/2} \sum_{k=1}^n \xi(k) ||G(\nu;r)
 \le C(\nu;r) ||\xi||G(\nu;r); \eqno(5.2)
$$

$$
\sup_i||\xi(i)||G\left(\nu_{\{\xi\}} \right) \le \sup_n ||n^{-1/2}
\sum_{k=1}^n \xi(k) ||G \left(\nu_{\{\xi \}} \right) \le
$$

$$
 C(\nu_{\xi}) \sup_i ||\xi(i)||G \left(\nu_{\{\xi \}} \right). \eqno(5.3)
$$
Recall that
$$
\nu_{\{\xi \}}(p) \stackrel{def}{=}  \sup_i ||\xi(i)||G(\nu;r).
$$

\vspace{3mm}
{\bf Proof.} Let $ \sup_i ||\xi(i)||G(\nu;r) < \infty $ and $ {\bf E}\xi = 0. $ We can
and will assume without loss of generality $ \sup_i ||\xi(i)||G(\nu;r)=1. $ It follows
from this equality

$$
\sup_i |\xi(i)|_p \le \nu(p), \ p \in [2,r).
$$
 We conclude by means of inequality (5.1)  by virtue of inequality $ p < r: $

$$
 \sup_n |n^{-1/2} \sum_{k=1}^n \xi(k) |_p  \le p \ \sqrt{2} \ |\xi|_p \le p \ \sqrt{2} \
 \nu(p) \le C_4(r) \nu(p), \eqno(5.4)
$$

$$
\sup_n \sup_{p \in [2, r)}  \left[ \left|n^{-1/2} \sum_{k=1}^n \xi(k) \right|_p /\nu(p) \right] \le C_4(r),
$$
which is equivalent to the assertion  of our theorem.\par
The second proposition of one follows from (5.4) after choosing $ \nu(p)= \sup_i \nu_{\xi(i)}(p). $\par

 As a consequence: \par

 \vspace{3mm}

 {\bf Theorem 5.2.}

 $$
 U(x) \le C(r,\Law( \{\xi(i) \})) \ T^{(\nu)}(x). \eqno(5.5)
 $$

\vspace{3mm}

{\bf Example 5.1.} Let $ {\bf E}\xi = 0 $ and suppose the r.v. $  \xi(i) $ satisfies the
condition (1.1) uniformly over $ i.$ Then for the values $ x > e $

 $$
 x^{-r} \ \log^{\gamma}(x) \ L(\log x) \le U(x) \le
 C_6(r,\gamma,L) \ x^{-r} \ \log^{\gamma+1}(x) \ L(\log x). \eqno(5.6)
 $$

\vspace{3mm}

\section{Superheavy tails}
\vspace{3mm}

We will investigate in this section the case when the r.v. $ \xi $ has symmetrical distribution
such that

$$
T_{\xi}(x) = \frac{K}{\log^{\kappa}x} L(\log x),  \eqno(6.1)
$$
where as before $ K, \kappa = \const > 0, L(z) $ is slowly varying as $ z \to \infty $ positive
continuous function, $ \xi(k), \ k = 1,2,\ldots $ are independent copies of $ \xi. $\par
{\it We do not know the exact norming sequence for the sums } $ \sum_{k=1}^n \xi(k). $  We
intend to represent here a more slight result as in the second section. \par
 Let us define for any fixed increasing deterministic tending to $ +\infty $ non-random
 sequence $ w(n), \ w(1) = 1 $ the following norming sequence $ B_n: $

$$
B_n = \exp \left( (K n)^{1/\kappa} \
L^{1/\kappa} \left( (K n)^{1/\kappa}  \right) \ w(n) \right) \eqno(6.2)
$$
and introduce correspondingly

$$
S(n) = \frac{\sum_{k=1}^n \xi(k)}{B_n},
$$

$$
\overline{U}(x) = \sup_n T_{S(n)}(x), \ x > e^2.
$$

\vspace{3mm}

{\bf Theorem 6.1.} For some positive finite constant $ C = C(\kappa, L, \{w(n)\}) $

$$
T_{\xi}(x) \le \overline{U}(x) \le T_{\xi}(x/C). \eqno(6.3a)
$$
 Further, the sequence $ S(n) $ tends in probability to zero as $ n \to \infty: $

 $$
\frac{\sum_{k=1}^n \xi(k)}{B_n} \stackrel{{\bf P}}{\to} 0, \eqno(6.3b)
 $$
 a "double Weak" Law of Large Numbers.\par

\vspace{3mm}

{\bf Proof.} We conclude as before  as $ t \to 0+ $
$$
\psi_{\xi}(t) = 2t \int_0^{\infty} \sin(tx) \ T_{\xi}(x) \ dx \sim
$$
$$
2 t \cdot K \int_e^{\infty} \sin(tx) \ \frac{K}{\log^{\kappa}x} \ L(\log x)  \ dx. \eqno(6.4)
$$
 The exact asymptotic of the last integral as $ t \to 0+ $ (Fourier transform) is
 calculated in the classical book of A.Zygmund \cite{Zygmund1}, p. 186-188:

 $$
 \psi_{\xi}(t) \sim K \ |\log t|^{-\kappa} \ L(|\log t|). \eqno(6.5)
 $$
 Calculating the characteristical function for the sequence $ S(n), $ we conclude
by means of estimate (2.2)
$$
\sup_n \left|\psi_{S(n)}(t) \right| \le K_1 \ |\log t|^{-\kappa} \ L(|\log t|), \ t \in (0,1/e),
\eqno(6.6)
$$
and

$$
\lim_{n \to \infty} \left|\psi_{S(n)}(t) \right| \to 0. \eqno(6.7)
$$
 The second proposition (6.3b) of theorem 6.1 follows immediately from (6.6),
 the first is proved alike the proof of theorem 2.1.\par

\vspace{3mm}

\section{Continuity and Stable Limit theorems for heavy tail random fields }
\vspace{3mm}

{\bf 1. Continuity.} \par
Let $ \eta(v), \ v \in V $  be  separable random field  (r.f.) (process)
defined aside from the probabilistic space $ \Omega $ on any set $  V. $ We suppose that
for arbitrary point $ v \in V $ the r.v. $ \eta(v) $ satisfies the condition (1.1) up to
continuous bilateral bounded multiplicative constant $ K(v): $

 $$
 T_{\eta(v)}(x) = K(v) \ x^{-r} \ \log^{\gamma}(x) \ L(\log x),  r = \const \in (1, \infty),
 \ x > e, \eqno(7.1)
 $$

$$
\gamma > -1, \ C_1 \le K(v) \le C_2, \ C_1, C_2 = \const, 0 \le C_1 \le C_2 < \infty.
$$
 Without loss of generality we can and do assume that for some fixed non-random
 value $ v_0, \ v_0 \in V \ K(v_0) = 1. $ \par
 Let us introduce the following function:

 $$
 \theta(p) = \left|\eta(v_0) \right|_p = \nu_{\eta(v_0)}(p), \ 1 \le p < r;
 $$
then

$$
\theta(p) \sim (r-p)^{-\gamma-1} \ L(1/(r-p)), \ p \to r-0.
$$

 From the equality (7.1) follows that

 $$
\sup_{v \in V} ||\eta(v)||G\theta  < \infty. \eqno(7.2)
 $$

 The so-called {\it natural distance} $ d(v_1,v_2) $ (more exactly, semi-distance:
 from the equality $ d(v_1,v_2)=0 $ does not follow $ v_1 = v_2) $ may be  defined
 by the formula

 $$
 d(v_1,v_2)= ||\eta(v_1) -  \eta(v_2)||G\theta. \eqno(7.3)
 $$
 The boundedness of $ d(v_1,v_2) $ follows immediately from (7.2). \par
{\bf  Remark 7.1.} The continuity of the coefficient $ K = K(v) $ is understood
relative the distance $ d = d(v_1,v_2). $ \par

 We denote as usually the metric entropy of the set $ V $ in the distance $ d(\cdot,\cdot) $
as $ H(V,d,\epsilon); $ recall that $ H(V,d,\epsilon) $ is the natural logarithm of the
minimal number of $ d - $ closed balls with radius $ \epsilon, \ \epsilon > 0 $ which
cover  the set $  V. $ By definition,

$$
N(V,d,\epsilon)= \exp H(V,d,\epsilon).
$$
The classical theorem of Hausdorff tell us that  $ \forall \epsilon > 0 \ N(V,d,\epsilon)< \infty $
iff the set $  V $ is precompact set relative the distance $  d. $ \par

\vspace{3mm}

{\bf Theorem 7.1.} If the following integral converges:

$$
\int_0^1 N^{1/r}(V,d,\epsilon)  \ H^{\gamma/r}(V,d,\epsilon) \ L^{1/r}(H(V,d,\epsilon))
\ d \epsilon < \infty, \eqno(7.4)
$$
then the trajectories $ \eta(t) $ are $ d(\cdot, \cdot) $ continuous with probability one:

$$
{\bf P}(\eta(\cdot) \in C(V,d)) = 1 \eqno(7.5)
$$
and moreover

$$
{\bf P}(\sup_{v \in V} |\eta(v)| \ge x) \le C \ \ x^{-r} \ \log^{\gamma}(x) \ L(\log x), \
x \ge e. \eqno(7.6)
$$

\vspace{3mm}

{\bf Proof.} From the definition of the norm $ ||\cdot||G\theta $ follows the inequality

$$
{\bf P} \left( \left| \frac{\eta(v_1) - \eta(v_2)}{d(v_1,v_2)} \right| \ge x \right) \le
C_1 \ x^{-r} \ \log^{\gamma}(x) \ L(\log x), \ x > e. \eqno(7.7)
$$
 By definition, we adopt in (7.7) that  $  0/0 = 0. $ \par
It remains to use  the result of the chapter 4, section 4.3 of the monograph
\cite{Ostrovsky1}. \par

\vspace{3mm}
{\bf  Remark 7.1.} The case when

$$
\sup_{v \in V} |\eta(v)|_r < \infty, r \ge 1
$$
and the distance

$$
\rho(v_1, v_2) = |\eta(v_1) - \eta(v_2) |_r
$$
was considered by G.Pizier \cite{Pizier1}. Indeed, if

$$
\int_0^1 N^{1/r}(V,\rho,\epsilon) \ d \epsilon < \infty,
$$
then

$$
{\bf P}(\eta(\cdot) \in C(V,\rho)) = 1
$$
and
$$
|\sup_{v \in V} |\eta(v)| \ |_r < \infty.
$$
{\bf Remark 7.1.} Another approach via the so-called majorizing measures, see in
\cite{Fernique1}, \cite{Talagrand1}, \cite{Talagrand2}.\par

\vspace{3mm}
{\bf 2. Stable  and Central Limit theorems.}\par
\vspace{3mm}

 We assume in addition that the random field $ \eta(v) $ is symmetrically distributed.
Let $ \eta_k(v) $ be independent copies of $ \eta(v). $ Define as before the norming
sequence $ b(n) $ as a solution of equation
$$
  n^{-1} = b^{-r}(n) \ |\log b(n)|^{\gamma} \ L(\log b(n))
$$
in the case $ r \le 2 $ and $ b(n) = \sqrt{n} $ when $ r > 2.$
 Put

 $$
 \beta_n(v)= \frac{1}{b(n)}\sum_{k=1}^n \eta_k(v). \eqno(7.8)
 $$

 The finite-dimensional distributions of the random fields $ \beta_n(v) $ converge
 as $ n \to \infty $ to the finite-dimensional distribution of the random field, which
 we denote as $ \beta(v). $ The last r.f. has a stable distribution when $ r < 2 $ and
 has a Gaussian centered distribution with at the same covariation function as $ \eta(v):$

$$
{\bf E} \beta(v_1) \beta(v_2) = {\bf E} \eta(v_1) \eta(v_2).
$$
 More information about stable distributions  in the Banach spaces
 see  in the monograph N.N.Vakhania, V.I.Tarieladze and S.A.Chobanan \cite{Vakhania1},
 chapter 5. \par
   We will say as ordinary that when the sequence r.f. $ \beta_n(\cdot) $  and r.f. $ \beta(\cdot)$
are $ d - $ continuous with probability one and the distributions in the space $ C(V,d) $
of r.f.  $ \beta_n(\cdot) $ converge weakly as $ n \to \infty $ to the distribution of
 $ \beta(\cdot), $ that the random fields $ \eta_k(v), \ k=1,2,\ldots $
 satisfy the Limit Theorem  in the space $ C(V,d). $ \par
 In the case $ r < 2 $ we have the Stable Limit Theorem; when $ r \ge 2 $ one can say as
 Central Limit Theorem in this space.  \par

  The limit theorems in the Banach spaces was initiated by Yu.V.Prokhorov
 in \cite{Prokhorov2} and was continued in many works, e.g.,  \cite{Billingsley1},
 \cite{Buldygin1}, \cite{Dudley1}, \cite{Mushtary1}. \par
 The applications of the limit theorems in Banach spaces in the Monte-Carlo method see in
\cite{Frolov1}, \cite{Grigorjeva1}.\par

\vspace{3mm}

{\bf Theorem 7.2.} If the following integral is finite:

$$
\int_0^1 N^{1/r}(V,d,\epsilon)  \ H^{(\gamma +1)/r}(V,d,\epsilon) \ L^{1/r}(H(V,d,\epsilon))
\ d \epsilon < \infty, \eqno(7.9)
$$
then the random fields $ \eta_k(v), \ k=1,2,\ldots $
satisfy the Limit Theorem  in the space $ C(V,d). $ \par
 Moreover,

$$
\sup_n {\bf P}(\sup_{v \in V} |\beta_n(v)| \ge x) \le C \ \ x^{-r} \ \log^{\gamma+1}(x) \ L(\log x), \
x \ge e. \eqno(7.10)
$$

\vspace{3mm}
{\bf Proof.} As long as the condition (7.10) is more strong as (7.4), we conclude
$ {\bf P}(\beta_n(\cdot) \in C(V,d)) = {\bf P}(\beta(\cdot) \in C(V,d)) = 1. $ \par
 It remains to prove the {\it tightness} of the measures $ \mu_n $ generated by the sequence
 $ \{ \beta_n(\cdot) \}: $

 $$
 \mu_n(A) = {\bf P} (\beta_n(\cdot) \in A),
 $$
 where $  A $ is Borelian set, in the space $ C(V,d). $ We obtain using theorem 2.1 that

$$
\sup_n{\bf P} \left( \left| \frac{\beta_n(v_1) - \beta_n(v_2)}{d(v_1,v_2)} \right| \ge x \right) \le
C_2 \ x^{-r} \ \log^{\gamma+1}(x) \ L(\log x), \ x > e. \eqno(7.11)
$$
 As before, we adopt by definition in (7.11) that $  0/0 = 0. $ \par
It remains to use  the result of the chapter 4, section 4.3 of the monograph
\cite{Ostrovsky1}. \par

\vspace{3mm}

{\bf 3. Applications.} \par

\vspace{3mm}
{\bf A.} We return now to the problem computation of (multiple) {\it parametric} integral of
a view:

$$
I(v) = \int_D f(v,y) \ \nu(dy), \ v \in V, \eqno(7.12)
$$
where $ \nu(\cdot) $ is again probabilistic measure on the set $ D: \ \nu(D) = 1.$ \par
 Let $ \tau(k), \ k=1,2,\ldots,n $ be as before independent r.v. with distribution
 $ \nu: \ {\bf P}(\tau(k) \in A) = \nu(A). $  The Monte-Carlo consistent estimation $ I_n(v) $
 of an integral  $ I(t) $ is

 $$
 I_n(v) = n^{-1} \sum_{k=1}^n f(v,\tau(k)). \eqno(7.13)
 $$

 Suppose for some $ r \in (1,2)$ and for all the values $ v \in V $
$$
 {\bf E} |f(v,\tau(1))|^r  < \infty
$$
or more generally that the r.v. $ f(v,\tau(k)) - I(v) $ satisfies the condition (1.1);
uniformly in $ v.$  In order to construct a non-asymptotical confidence interval for $ I $ of
a reliability $ 1 - \delta, \ \delta = 0.05; 0.01 $ etc. in uniform over $ v \in V $
we consider the probability

 $$
 U_n(x) = {\bf P} \left( \sup_{v \in V} b(n)^{-1} \left|\sum_{k=1}^n (f(v,\tau(k)-I(v))) \right| > x \right).
 \eqno(7.14)
 $$

Note that if the sequence of r.f. $ \{ f(v,\tau(k)-I(v)) \} $ satisfy Limit Theorem, then

$$
\forall x > 0 \ \lim_{n \to \infty} U_n(x) \to U(x), \eqno(7.15)
$$
where

$$
 U(x) = {\bf P} \left( \sup_{v \in V} |\zeta(v)| > x \right), \eqno(7.16)
 $$
$ \zeta(v) = \zeta(\omega,v) $ is stable or Gaussian random field.\par
The asymptotical or non-asymptotical behavior of $ U(x) $ as $ x \to \infty $
in both the cases: SLT or CLT is known, see, e.g. \cite{Vakhania1}, chapter 5; 
\cite{Piterbarg1}, \cite{Ostrovsky1}, chapter 3.\par

Therefore, we conclude asymptotically as $ n \to \infty $ denoting by $ X(\delta) $ the 
solution of an equation

$$
U(X(\delta)) = \delta
$$
that with probability at least $ 1-\delta $ in the uniform norm

$$
\sup_{v \in V} |I_n(v) - I(v)| \le X(\delta) b(n)/n. \eqno(7.17)
$$
 Note that as before
 $$
  \lim_{n \to \infty} b(n)/n = 0.
 $$

\vspace{3mm}

 Notice that it may be used {\it non-asymptotical approach,} where the probability $ U_n(x) $
 allows the evaluating as follows:

 $$
 U_n(x) \le \sup_n U_n(x) \le C \ x^{-r} \ \log^{\gamma+1}(x) \ L(\log x). \eqno(7.18)
 $$

\vspace{3mm}

{\bf  B.} Analogous application appears in statistics. Indeed, let us consider the following
classical  scheme of date-process elaboration.

$$
\eta_k(v) = \theta(v) + \xi_k(v), \ k=1,2,\ldots,n; \eqno(7.19)
$$
where $ \theta(v), \ v \in V $ is unknown deterministic function, $ \{ \xi(k) \} $ are i.,i.d.
centered r.f. satisfying the condition (1.1) with $ r > 1 $ (additive noise with heavy
tail). \par
 The consistent estimation of the functional parameter $ \theta(v) $ has a view

 $$
\hat{\theta}_n(v) = n^{-1} \sum_{k=1}^n \eta_k(v). \eqno(7.20)
 $$
We conclude as before that with probability at least $ 1-\delta $
under formulated above conditions and notations  

$$
\sup_{v \in V}\left|\hat{\theta}_n(v) - \theta(v) \right| \le X(\delta) b(n)/n. \eqno(7.21)
$$

\vspace{3mm}

\section{Concluding remarks}
\vspace{3mm}

{\bf A. Non-symmetrical case.}\par

\vspace{4mm}
 The results of the second and third section remains true still without restriction of
symmetrical distribution of the independent r.v. $ \{\xi(i) \}. $ Indeed, it is sufficient
to assume in the case $ r \in (1,2] $ in addition to the equality (1.1)
$$
{\bf E} \xi(i)=0
$$
and in the case $ r=1 $
$$
\sup_{a > 0} {\bf E}|\xi(i) \ I(|\xi(i)| \le a) | < \infty;
$$
see, e.g., \cite{Bahr1}, \cite{Braverman2}. The proof may be obtained also from the
symmetrization arguments, see  also \cite{Latala1}. \par
\vspace{3mm}

{\bf B. Not identical distributed r.v.} \par

\vspace{3mm}

It is not necessary to suppose also in the independent case when $ r > 2 $ as in the
remark 4.2  that the r.v. $ \xi(i) $ are identical distributed; it is sufficient to
assume in addition that the r.v. $ \{\xi(i) \} $ are independent, centered and  such that
for some positive finite  constants $ C_1, \ C_2 $
$$
 C_1 \ x^{-r} \ \log^{\gamma}(x) \ L(\log x) \le
 T_{\xi(i)}(x) \le C_2 \ x^{-r} \ \log^{\gamma}(x) \ L(\log x), \ x > e.
 $$

\vspace{3mm}
{\bf C. About calculation of the norming sequence.}\par
\vspace{3mm}
 We investigate here the equation $ n \psi(1/b(n)) = 1 $ for the norming sequence
 $ \{ b(n)  \} $ in the case $ r < 2. $ \par
  We have under condition (1.1):

  $$
  n^{-1} = b^{-r}(n) \ |\log b(n)|^{\gamma} \ L(\log b(n)). \eqno(8.1)
  $$
 It is reasonable to assume that  as $ n \to \infty $

$$
b(n) \sim n^{1/r} \ \log^{\gamma/r}(n) \ L^{1/r}(\log n). \eqno(8.2)
$$
 We obtain substituting into (6.1C) analogously to the classical monograph of E. Seneta
 \cite{Seneta1}, p. 29-32 that the asymptotical expression for $ b(n) $ is true when the
 following condition holds:

 $$
 L \left( \frac{X^{1/r}}{\log^{\gamma/r} X \cdot L^{1/r}(\log X) } \right)
 \asymp L(X), \ X \to \infty. \eqno(8.3)
 $$
 Note that the condition 6.3C is satisfied for the function, e.g.,

 $$
 L(X) = C \ (\log X)^{\Delta},  \ \Delta = \const.
 $$

\vspace{3mm}

{\bf D. Tail comparison through moments inequalities.}\par

\vspace{3mm}

 If for two r.v. $ \xi $ and $ \eta \  T_{\xi}(x) \le T_{\eta}(x), \ x > 0, $
then evidently
$$
|\xi|_p \le |\eta|_p, \ p \ge 0. \eqno(8.4)
$$
 We discuss in this pilcrow the inverse problem. Indeed, suppose the inequality (6.1) holds.
Our purpose is to obtain the estimate the upper bound for tail probability $ T_{\xi}(x). $ \par
 In detail, assume that for {\it any values} $ p $ from the non-trivial segment $ p \in [1,r),
  \ r \in (1,\infty) $

$$
|\xi|_p \le |\eta|_p, \ p \ge 0. \eqno(8.5)
$$

 It follows from Tchebychev's inequality

 $$
 T_{\xi}(x) \le x^{-p} \ |\eta|_p^p, \ x > 0, \ p \in \supp \nu_{\eta},
 $$
 therefore

 $$
 T_{\xi}(x) \le \inf_{p \in \supp \nu_{\eta} }  \left[ x^{-p} \ |\eta|_p^p \right], \ x > 0. \eqno(8.6)
 $$
  If for instance $ \supp \nu_{\eta} = [1,r),  $ or equally when $ \nu_{\eta}(r+0) = \infty $
  for some $ r > 1, $  then

 $$
 T_{\xi}(x) \le \inf_{p \in [1,r) }   \left[ x^{-p} \ |\eta|_p^p \right], \ x > 0. \eqno(8.7)
 $$
 We get, e.g. in particular choosing in (9.7) the value $ p = r - C/\log x, \ C = \const > 0 $
 for sufficiently greatest values $ x $

 $$
 T_{\xi}(x) \le e^C \ x^{-r} \ {\bf E} |\eta|^{r-C/\log x}.  \eqno(8.8)
 $$

\vspace{3mm}

{\bf E.  Non-uniform norming sequence.} \par

\vspace{3mm}

 M.Braverman in the article \cite{Braverman2} considered a more general as uniform
norming vector $ a = \{ a(1), a(2), \ldots, a(n)\}. $  In detail, let $ r \in (0,2) $
and let $ \{ \xi(k) \}, \ k=1,2,\ldots,n $ be again the independent copies of the
symmetrical r.v. $ \xi $ satisfying the condition (1.1). Put

$$
U^{(a)}(x) =  {\bf P} \left( \left|\sum_{k=1}^n a(k)\xi(k) \right| \right). \eqno(8.9)
$$

M.Braverman in \cite{Braverman2} introduced an Orlicz space of numerical sequences
$ a = \{ a(1), a(2), \ldots, a(n)\}  $ by means of Orlicz's function $ \psi(t) =
\psi_{\xi}(t): $

$$
||a||_{\psi} = \inf \{t, t > 0, \sum_{k=1}^n \psi(|a(k)|/t) \le 1 \}, \eqno(8.10)
$$
and proved that the unit ball in this space is natural norming sequence in the
sense  of $ L_p, \ p < r $ boundary.  \par

 Denote

 $$
 \overline{U}(x) = \sup_{ a: ||a||_{\psi} \le 1 } U^{(a)}(x). \eqno(8.11)
 $$
It may be proved analogously theorem 2.1 that

$$
\overline{U}(x) \le  0.5  x \int_{-2/x}^{2/x} \overline{\psi}(t) \ dt, \ x > e. \eqno(8.12)
$$

\end{document}